\documentclass[10 pt]{amsart}
\usepackage{amscd,amsmath,amsthm,amssymb}
\usepackage{enumerate,varioref}
\usepackage{epsfig}
\usepackage{graphicx}
\newtheorem{thm}{Theorem}

\newtheorem{prop}[thm]{Proposition}
\theoremstyle{definition}

\theoremstyle{remark}

\newtheorem{rem}[thm]{Remark}

\numberwithin{equation}{subsection}
\newcommand{\R}{\mathbb{R}}

\newcommand{\C}{\mathbb{C}}

\begin{document}

\title{The First Order Effect of the Quantum Weyl Algebra on a Harmonic Oscillator }
\author{Clark Alexander}

\email{clark@imsc.res.in} \maketitle \tableofcontents
\section{Introduction}

The goal of the present article is to give another example than those first presented in a pair of papers from the mid 1990s.  Since then these papers, [CJSS] and [DaP], have received little attention.  Both papers explored the idea of a non-relativisitic free particle in a noncommutative space.  The former explored only a deformed space, whilst the latter explored a deformed space time.  In this article we shall consider only a deformed space.  Moreover, we shall consider only a deformed three-dimensional space.  In what follows, we shall explore the effect of a noncommutative space on the quantum harmonic oscillator.  The deformation present in this article has partly fallen out of favor because it fails to obey a $SO_q(3)$ symmetry and in the case of [DaP] it fails to obey $SO_q(3,1)$ symmetry.  However, to begin we shall show that while quantum orthogonal symmetries are not obeyed, there is a different symmetry, that of $Sp_q(6)$ which the Quantum Weyl Algebra does obey.

\section{Preliminaries on the Quantum Weyl Algebra}

The quantum Weyl algebra as it is now so called was first seen explicitly in [PW].  It is now given as a covariant differential calculus on $\mathcal{O}(\C^N_q)$ named $\mathcal{A}_q(N)$ generated
by $\{X_j,\partial_j | i,j=1,\dots N\}$  with the following relations:

\begin{eqnarray}
X_i X_j &=& q X_j X_i,\hspace{5mm}  i<j,\\
\partial_i \partial_j &=& q^{-1} \partial_j \partial_i,\hspace{5mm}  i<j, \nonumber\\
\partial_i X_j &=& q X_j \partial_i,\hspace{5mm}  i\neq j, \nonumber\\
\partial_i X_i - q^2 X_i \partial_i &=& 1 + (q^2-1)\sum_{j>i}X_j\partial_j. \nonumber 
\end{eqnarray}

\begin{rem}
This covariant differential calculus was also seen in [WZ], however this was given only in the case $N=2$.
\end{rem}

\subsection{Quantum Symmetries}

In both [CJSS] and [DaP] similar comments are made on the fact that $\mathcal{A}_q(3)$ does not obey $SO_q(3)$ symmetry.  It requires very little to see that these statements are indeed valid.  Fortunately, the subsets 
$\{ X_i | i=1,\dots, N\}$ and $\{ \partial_i | i=1,\dots,N\}$ obey $GL_q(N)$ and $GL_{q^{-1}}(N)$ symmetries respectively.  If one dives into the literature ever so slightly deeper, one can find the presentation of [KS] in \S 12.3.3 exhibiting 
$\mathcal{A}_q(N)$ as a left quantum space of $Sp_q(2N)$.  

For much of this article we shall only be concerned with $N=3$ and so we will now give a brief synopsis of how $Sp_q(6)$ acts on $\mathcal{A}_q(3)$.

Let us redefine the elements of $\mathcal{A}_q(3)$ as $\{ y_i | i=1,\dots,6\}$ where 
\begin{eqnarray}
y_1 = \alpha q \partial_3, & y_2 = \alpha q^2 \partial_2, & y_3 = \alpha q^3 \partial_1, \\
y_4 = X_1, & y_5 = X_2, & y_6 = X_3, \nonumber
\end{eqnarray}
where $\alpha \in \C^{\times}$.

\begin{prop}
Let $\hat R = (\hat R^{ij}_{kl})$ be the $R$-matrix and $C=(C^i_j)$ the metric for $Sp_q(6)$.  Then the generators $\{X_1,X_2,X_3,\partial_1,\partial_2,\partial_3\}$ satisfy 
the relations for $\mathcal{A}_q(3)$ if and only if
\begin{equation}
\sum_{k,l}^6 \hat R^{ij}_{kl}y_ky_l - q y_iy_j - \alpha q^{-3}C^i_j = 0, \hspace{5mm} i,j=1,\dots,6.
\end{equation}
\end{prop}

The $R$-matrix and metric are given explicitly in [KS] \S 9.3.
The proof of this proposition is essentially a direct calculation and verification that the above  equation reduces to
\begin{equation}
y_{4-j}y_j - q^{-2}y_j y_{4-j} = -q^{-j}\alpha (q^{-2}-1) \sum_{k<j}q^{k-j}y_k y_{4-k}.
\end{equation}

\subsection{A Concrete Realization of $\mathcal{A}_q(3)$}
Since we shall only be working within a three dimensional system, for the remainder of this paper we will adopt the notations that capital letters $(X,Y,Z)$ will represent noncommutative coordinates and lower case letters $(x,y,z)$ will represent commutative coordinates.  Furthermore we will also use $X_1 = X, X_2=Y,X_3=Z$ as well as $x_1=x,x_2=y,x_3=z$ and use the two interchangeably when convenient.

Define the following quantities
\begin{eqnarray}
q &:=& e^{i\theta}, \\
M_j &: =&  x_j \partial_{x_j}, \nonumber\\
M_{k>j} &:=& \sum_{k>j}M_j,\nonumber\\
\beta_j &:=& \left\{ \frac{q^{2(M_j+1)}-1}{(q^2-1)(M_j+1)} \right\}^{1/2}.\nonumber
\end{eqnarray}
Here we restrict to $\theta\in\R$.  We also note that in [CJSS] and [DaP] the operators $M_j$ are written as $N_j$.  However in this paper, as a manner of avoiding confusion, 
we shall use $M_j$ since $N_j$ generally represent number operators in quantum mechanics.

In this scenario we now may write our noncommutative coordinate system explicitly as

\begin{eqnarray}
X_j & = & x_j \beta_j  q^{M_{k>j}},\\
\partial_{X_j} &=& q^{M_{k>j}} \beta_j \partial_{x_j}.
\end{eqnarray}

Expanding this for the sake of clarity we see
\begin{eqnarray}
X = x \beta_1 q^{M_2+M_3},                              &   Y = y \beta_2 q^{M_3},                              & Z = z\beta_3,\\
\partial_X = q^{M_2+M_3} \beta_1 \partial_x, &  \partial_Y = q^{M_3} \beta_2 \partial_y,  & \partial_Z = \beta_3 \partial_z.
\end{eqnarray}

It is a relatively easy computation to verify that these operators satisfy the relations of our algebra $\mathcal{A}_q(3)$.

\section{Slightly Noncommutative Space}

Looking back on [CJSS] and [DaP] we have a definiton of a \emph{slightly noncommutative space} (SNCS).  In the case of [DaP] they deform time as well for a SNCST.  This space is given
by assuming that the deformation parameter $q$ is close to 1, or more importantly, assuming $\theta \approx 0$.  With this in mind, our parameter becomes
\[
q \approx 1 + i\theta.
\]

The consequences of this are mildly startling.  Of course, our deformation space and commutation relations change drastically, but we may now write down coordinates in a new explicit form.
Let us begin with $\beta_j$.  In the first order, the expansion is rather easy, not offering such complications as in the second order expansion in [A].  

In order to expand $\beta_j$ properly it is best to consider the denominator first.  In the denominator, allowing a term of order $\theta$ equates to a term of order $\theta^{-1}$ overall.  For this
reason it will be important to carry a term of order $\theta^2$ in the numerator.
Therefore the denominator gives
\[
(q^2-1)(M_j+1) \approx (1 + 2i\theta -1)(M_j+1) = (2i\theta)(M_j+1),
\]
and the numerator yields
\begin{eqnarray*}
q^{2(M_j+1)}-1 & \approx &  1 + 2i\theta(M_j+1) + \frac{1}{2}(2i\theta(M_j+1))^2 - 1\\ 
                            &     =        &2i\theta(M_j+1)(1 + i\theta(M_j+1)).
\end{eqnarray*}

Putting this together with the expansion 
\[
\sqrt{1+x} \approx 1 + \frac{1}{2}x + O(x^2)
\]
we recover

\begin{equation}
\beta_j = \left\{ \frac{q^{2(M_j+1)}-1}{(q^2-1)(M_j+1)} \right\} ^{1/2} \approx 1 + \frac{1}{2}i\theta (M_j+1).
\end{equation}

Keeping track of $q^{M_{k>j}} \approx 1 + i\theta M_{k>j} $ we have the following relations up to first order:

\begin{eqnarray}
\partial_X & \approx & [1+ \frac{1}{2} i\theta (M_1+1) + i\theta (M_2+M_3)]\partial_x,\\
\partial_Y & \approx & [1+ \frac{1}{2} i\theta (M_2+1) + i\theta M_3] \partial_y, \nonumber \\
\partial_Z & \approx & [1+ \frac{1}{2} i\theta (M_3+1)] \partial_z, \nonumber \\
X              &  \approx & x [1+ \frac{1}{2} i\theta (M_1+1) + i\theta (M_2+M_3)], \nonumber \\
Y              & \approx & y[1+ \frac{1}{2} i\theta (M_2+1) + i\theta M_3], \nonumber \\
Z              & \approx & z[1+ \frac{1}{2} i\theta (M_3+1)]. \nonumber
\end{eqnarray} 

\section{The Effect on the Ground State of the Quantum Harmonic Oscillator}

Up until now we have not explored the idea of making a quantum mechanical system out of our algebra or out SNCS.  The procedure adopted by both [CJSS] and [DaP] is to simply
map a Hamiltonian in commutative coordinates into a Hamiltonian in noncommutative coordinates by a simple procedure which looks as follows:

\begin{eqnarray}
x_j &\mapsto & X_j,\\
\partial_{x_j}&\mapsto & \partial_{X_j},\nonumber \\
P_{X_j} & = & -i\hbar \partial_{X_j}, \nonumber\\
H(p_{x_j},x_j) & \mapsto & \mathcal{H}(P_{X_j},X_j).\nonumber
\end{eqnarray}

Given this procedure, we now write the noncommutative Hamiltonian of a quantum harmonic oscillator as
\begin{equation}
\mathcal{H} = \frac{1}{2m}(P_X^2 + P_Y^2 + P_Z^2) + \frac{1}{2}\omega^2 (X^2 +Y^2 + Z^2).
\end{equation}

Here, we will allow ourselves the mathematical luxuries of assuming a unit mass, unit frequency, and unit Planck constant ($m=1,\omega=1, \hbar =1$).

In this new context, the equation we wish to solve is 
\begin{equation}
i \frac{\partial \Psi}{\partial t} = \frac{1}{2}(P_X^2+P_Y^2+P_Z^2 +X^2 + Y^2 + Z^2)\Psi.
\end{equation}

It is not known whether or not this system has a general solution which is unique.  The most natural thing one can attempt 
in this scenario is successive approximations with higher order terms.  The mathematics which follows requires some trickery and 
is at first seemingly unnatural.  In the calculation of the effect of this system on a free particle [CJSS] and [DaP] assume the solution of a slowly varying free particle.
Here, we shall mimic this by assuming the solution of a harmonic oscillator in its ground state.

That is
\begin{equation}
\Psi = \pi^{-1/4} e^{-r^2/2},
\end{equation}

where $\pi^{-1/4}$ is a normalization constant and $r=(x,y,z)$ is the radial vector.  The intuition behind this assumed solution is that thie $\Psi$ is the time independent solution to
\[
\frac{1}{2}(p_x^2+p_y^2+p_z^2+x^2+y^2+z^2)\Psi = E \Psi.
\]

Since our coordinate system is slightly noncommutative the assumption we make is that we begin in the standard quantum mechanical system and then slowly ``turn on" the noncommutative coordinate system.
The question we are asking is, ``what happens to the physical system as this adiabatic change happens?"

To answer this question we will proceed in the following way:  Expand the noncommutative operators to first order, let them hit the assumed solution, and collect the ``new" terms.  That is to say, in our expansion of 
$X$ we have the form $x (1 + \theta \bullet)$ wherein we will leave the term $x$ alone and collect the other terms.  This will allow us to see the changes in the system.  Given this, we obtain

\begin{eqnarray} 
\partial_X \Psi & \approx & \partial_x \Psi + i\theta x (\frac{1}{2}x^2 + y^2+z^2 -1)\Psi, \\
\partial_Y \Psi & \approx & \partial_y \Psi + i\theta y (\frac{1}{2}y^2+ z^2 -1)\Psi, \nonumber\\
\partial_Z \Psi & \approx & \partial_z \Psi + i\theta z (\frac{1}{2}z^2 -1)\Psi, \nonumber \\
              X \Psi & \approx & x \Psi - i\theta x (\frac{1}{2}(x^2-1) + y^2+z^2)\Psi, \nonumber\\
              Y \Psi & \approx & y \Psi -i\theta y (\frac{1}{2}(y^2-1) + z^2) \Psi, \nonumber \\
              Z \Psi & \approx & z \Psi - i\theta z (\frac{1}{2}(z^2-1)) \Psi. \nonumber 
\end{eqnarray}

From here, we simply let $P_{X_j} = -i\partial_{X_j}$ and find that our Hamiltonian has taken the form

\begin{equation}
\mathcal{H}\Psi \approx (\frac{1}{2} \sum_{j=1}^3 (p_{x_j} - A_{x_j})^2 + V_{R} + i V_{I})\Psi.
\end{equation}

We have what looks to be a magnetic potential as well as a dissipative term.  Here, the explicit formulae are

\begin{eqnarray}
A_x & = & \theta x (\frac{1}{2}x^2+y^2+z^2-1), \\
A_y & = & \theta y(\frac{1}{2}y^2 + z^2 -1), \nonumber\\
A_z & = & \theta z(\frac{1}{2}z^2-1), \nonumber \\
V_{I}& = & -\theta[(x^3+y^3+z^3) + (xy^2 + xz^2 + yz^2) + (x+y+z)]. \nonumber
\end{eqnarray}

yielding a magnetic field of 
\begin{eqnarray}
B_x & = & -2\theta yz, \\
B_y & = & 2\theta xz, \nonumber \\
B_z & = & -2\theta yz. \nonumber
\end{eqnarray}

More succinctly,
\begin{equation}
B_i = \epsilon_{ijk}2\theta x_j x_k.
\end{equation}

\section*{Discussion}
One notes the striking resemblance to the results found in [CJSS] and [DaP] with respect to the form of the magnetic field.  However, in both of those works, an additional imaginary term appears in the momenta operators.
The imaginary component in the momenta disappears when considering the second order correction (cf [A]), but in the case of the harmonic oscillator, no such imaginary component appears.  Furthermore, another
interesting physical phenomenon occurs with the imaginary part of the potential.  In this case, we recover an "unstable" harmonic oscillator.  In a standard quantum mechanical setting (i.e. a textbook problem) one will find that a 
constant imaginary term added to the potential in the form $V = V_0 - i\alpha$ for $\alpha >0$ the probability of ``finding" a particle decreases in time.  Moreoever, one can easily deduce that the following:

\begin{eqnarray*}
H & = & T + V - i\alpha,\\
P(t) & = & \int |\Psi|^2 dt,\\
\frac{dP(t)}{dt} &=& \frac{-2\alpha}{\hbar}P(t).
\end{eqnarray*}
 
 In our case, we have a position dependent imaginary potential which serves to increase the probability in certain places, and decrease in others.  What we have left out in this article are the calculations and results for
 excited harmonic oscillators.  What one should find after performing these calculations is that a harmonic oscillator in state $|n_1,n_2,n_3\rangle$ does not stay in this state, but rather a mixing phenomenon occurs.
 In particular we should expect something of the form
 \[
 \mathcal{H}|n_1,n_2,n_3\rangle \approx \sum_{\alpha_i\in \{-1,0,1\} } E |n_1+\alpha_1,n_2+\alpha_2,n_3+\alpha_3\rangle.
 \]

\section*{Acknowledgements}
The author wishes to thank Professor R. Jagannathan for his suggestion of tackling this question.   Thanks are also in order for E. Rogers for his help in reading the preliminary work of this article and helping to explain when the physical interpretation was incorrect.

\end{document}